\documentclass[a4]{article}

\usepackage{amsfonts}
\usepackage{amsmath}
\usepackage{bm}
\usepackage{bbm}
\usepackage{graphicx}
\usepackage{nameref,hyperref}

\newcommand{\be}{\boldsymbol{e}}

\newcommand{\bl}{\boldsymbol{l}}

\newcommand{\bu}{\boldsymbol{u}}
\newcommand{\bx}{\boldsymbol{x}}
\newcommand{\by}{\boldsymbol{y}}
\newcommand{\bz}{\boldsymbol{z}}

\newcommand{\bX}{\boldsymbol{X}}

\newcommand{\bdelta}{\boldsymbol{\delta}}

\newcommand{\bphi}{\boldsymbol{\phi}}
\newcommand{\blambda}{\boldsymbol{\lambda}}
\newcommand{\bmu}{\boldsymbol{\mu}}
\newcommand{\btheta}{\boldsymbol{\theta}}

\newcommand{\R}{{\mathbb R}}
\newcommand{\Z}{{\mathbb Z}}

\begin{document}
\title{Computing the extinction path for epidemic models}
\author{Damian Clancy and John J.~H.~Stewart\\
Department of Actuarial Mathematics and Statistics\\ Maxwell Institute for Mathematical Sciences\\ Heriot-Watt University\\ Edinburgh\\ EH14 4AS\\ UK\\
d.clancy@hw.ac.uk}
\maketitle

\section*{Abstract}
In infectious disease modelling, the expected time from endemicity to extinction (of infection) may be analysed via WKB approximation, a method with origins in mathematical physics.
The method is very general, but its uptake to date may have been limited by the practical difficulties of implementation.
It is necessary to compute a trajectory of a (high dimensional) dynamical system, the `extinction path', and this trajectory is maximally sensitive to small perturbations, making numerical computation challenging.
Our objective here is to make this methodology more accessible by presenting four computational algorithms, with associated Matlab code, together with discussion of various ways in which the algorithms may be tuned to achieve satisfactory convergence.
We illustrate our methods using three standard infectious disease models.
For each such model, we demonstrate that our algorithms are able to improve upon previously available results.

\section*{Introduction}
A quantity of great interest in epidemiology is the expected time until infection dies out from a population.
In some cases, global eradication may be a target---as of 2024, the International Task Force for Disease Eradication~\cite{ITFDE} lists eight diseases that could potentially be eradicated globally: Guinea worm (dracunculiasis), poliomyelitis, mumps, rubella, lymphatic filariasis, cysticercosis, measles, and yaws.
Alternatively, interest may be in elimination from some local region, as during the 2001 outbreak of foot and mouth disease in the UK~\cite{KWS01}.
Where long-term elimination is not seen as feasible, the duration of each individual outbreak is of interest, as with Ebola virus disease in West Africa~\cite{AMF16,NBKFB17}, or plague in Madagascar~\cite{NPH18}.

If the basic reproduction number $R_0$ (the expected number
of secondary cases directly generated by a typical primary case in an otherwise susceptible population) is less than~1, then only minor outbreaks are possible (the process is subcritical).
The entire course of the infection process may then be modelled using a (linear) branching process~\cite{B83}, and the distribution of extinction time approximated using the approximating branching process~\cite{CT18}.
In the supercritical case $R_0 > 1$, following invasion into a na\"ive population, the infection may become endemic in the population, at which point the branching process no longer provides an appropriate model.
Random fluctuations around the endemic level may then lead to eventual extinction of infection.

A general approach to studying the expected extinction time from endemicity is via the WKB (Wentzel, Kramers, Brillouin) method.
This approach has its origins in mathematical physics (see, for example, chapter~10 of~\cite{BO99}), and numerous applications of the method to epidemic models have appeared in the mathematical physics literature, see~\cite{AM10,AM17} and references therein.
To date, there has been relatively little uptake of this methodology within the broader epidemic modelling and mathematical biology communities, with a few exceptions, e.g.~\cite{OM10,NBKFB17,CS24}.
Part of the explanation for this may be the difficulties in implementing the method in practice.
It is necessary to compute a trajectory of a (high dimensional) dynamical system specific to the epidemic model of interest, the `extinction path'.
The extinction path is defined over an infinite time interval, and is maximally sensitive to small perturbations~\cite{FBSS11,SFBS11}.
Considerable effort is therefore required in `tuning' computer code to obtain satisfactory convergence.
This is well illustrated in~\cite{BFB16}, where the WKB approach is successfully applied to a number of epidemic models, but substantial careful thought is required in each individual case.

The aim of this paper is take a step towards making WKB methodology accessible to applied researchers.
To this end, we present two basic approaches to computing the extinction path: (i)~a finite difference method; and (ii)~a collocation method.
For each of these methods, we consider two approaches to dealing with the infinite time interval over which the extinction path is defined: (i)~truncation to a finite time interval; and (ii)~transformation to a finite time interval.
We thus consider a total of four algorithms.
A finite difference method with time truncation was proposed in~\cite{LS13}, and has become the standard approach within the literature~\cite{BFB16,HS16,NBKFB17,HSS18,KHA22}.
Collocation and time transformation have been applied together in~\cite{CS24}, and time transformation is mentioned in~\cite{BFB16}, but so far as we are aware, these are the only instances to date of the application of either collocation methods or time transformations in computing epidemic extinction paths.

We illustrate our methods in a number of specific applications, and discuss a variety of ways in which the algorithms may be tuned to suit the model under consideration.
The illustrative applications that we consider are (i)~a stochastic Ross-Macdonald malaria model; (ii)~a network susceptible-infectious-susceptible (SIS) model; and (iii)~a susceptible-exposed-infectious-removed (SEIR) model.
We present Matlab code that may be modified to apply to any epidemic model of interest.
For the finite difference method with time truncation, Matlab code has previously been made available by~\cite{HSS18}; our code builds upon ideas present in the code of~\cite{HSS18}, with additional enhancements.
For the finite difference method with time transformation, and for the collocation method with truncation or transformation, the code presented here is, so far as we are aware, the first to be made available.
While it is straightforward in principle to modify our code to apply to other epidemic models, there may be considerable work required in tuning the code to suit the model; the discussion around our illustrative applications can provide guidance here.
For all of our illustrative applications, we present results that go well beyond previously available results for these models.

We find that all four of our algorithms can work well, even in high dimensions, provided the extinction path is of a simple shape.
Convergence becomes more difficult to achieve, requiring careful adjustment of tuning parameters, as the basic reproduction number $R_0$ increases far above~1; as the dimensionality of the problem increases; and when the extinction path is of more complicated shape.
It then becomes very useful to have four algorithms available, since any one of the algorithms may prove more effective than the others for a particular model of interest.

In the epidemic modelling context, it is usually desirable to estimate the expected extinction time across a range of parameter values, to allow for uncertainty as to true parameter values, and to model the effects of interventions.
Previous authors have often presented results from WKB methodology only for one set of parameter values, e.g.~\cite{NBKFB17}.
We present, for each of our illustrative examples, results across ranges of parameter values.
In the case of the finite difference method, computations across a range of parameter values are greatly facilitated through `vectorization' of our Matlab code~\cite{vectorization}, which results in considerably faster execution times than non-vectorized code such as that presented in~\cite{HSS18}.

All of our algorithms require the numerical solver to be supplied with an `initial guess' for the extinction path.
The standard way to do this involves parameter continuation~\cite{SKR00}, which we discuss in our~\nameref{methods} section below, and apply to the Ross-Macdonald malaria model and the SEIR model.
For the network SIS model, we instead propose a new approach: to generate our initial guess from the solution to a closely related, analytically solvable, problem.
This approach does not seem to have been used before; the recent results of~\cite{C24} open up the possibility of its more widespread use in the future.

The rest of the paper is organized as follows.
In the~\nameref{methods} section, we set out the steps in the WKB approach to estimating extinction times, survey various methods that have been proposed to solve the resulting Hamilton-Jacobi partial differential equation, and discuss some issues of implementation.
Next, in the~\nameref{models} section, we describe our three example models.
The~\nameref{results} section presents numerical results, obtained using our four algorithms, for each of our three example models, together with discussion of how the algorithms may be tuned to obtain satisfactory convergence.
Finally, in the~\nameref{discussion} section, we summarise our main results, discuss the extent to which we have succeeded in our aim of making WKB methodology more readily accessible, and suggest some possible directions for further work.
The Matlab code used to generate our results is presented in~\nameref{S1_File},~\nameref{S2_File},~\nameref{S3_File}.

\section*{Methods}\label{methods}
\subsection*{General theory}
We briefly set out the steps in the WKB approach, as it applies to the analysis of expected extinction time for epidemic models.
The method is well established in the literature, including comparisons with other approaches and with results from Monte Carlo simulation---see, for example,~\cite{CT18,AM10,AM17,CS24}.
A sketch justification is given in~\nameref{S1_Appendix}; for further details and more complete justification, see~\cite{AM10,AM17,C18a,C24}.

We suppose that the epidemic process is modelled  as a continuous-time Markov process on $\Z_+^k$ whose components represent numbers of different types of individuals (susceptible, infectious, etc.).
Suppose that the state space~$S \subseteq \Z_+^k$ may be partitioned as $S = C \cup \bar C$, where $\bar C$ consists of the disease-free states.
We assume that $C$ is a transient communicating class, and that the process will hit $\bar C$ within finite time with probability~1.
In general, $\bar C$ is of the form $\bar C = \{ \bx \in \Z_+^k : x_i = 0 \mbox{ for all } i \in \mathcal{I} \}$, where $\mathcal{I} \subseteq \{ 1,2,\ldots,k \}$ is the set of components corresponding to (different types of) infected individuals.

To obtain an approximation valid for large populations, we consider a sequence of Markov processes $\{ \bX^{(N)} (t) : t \ge 0 \}$ indexed by~$N$, where $N$~represents `typical' population size, and assume that this sequence of processes is density dependent in the sense of chapter~11 of~\cite{EK05}.
That is, transition rates are of the form
\begin{eqnarray}
\Pr \left( \bX^{(N)} ( t + \delta t ) = \bx + \bl \mid \bX^{(N)} (t) = \bx \right) 
&=& 
N \beta_{\bl} \left( \frac{\bx}{N} \right) + o ( \delta t )
\mbox{ as } \delta t \to 0
\label{density_dependent_rates}
\end{eqnarray}
for $\bx \in S$, $\bl \in \mathcal{L}$, where $\mathcal{L}$ is a finite set consisting of the possible jumps from each state $\bx \in S$, and for each $\bl \in \mathcal{L}$, $\beta_{\bl} ( \cdot )$ is a function from $\R_+^k$ to $\R_+$.
Then if $\lim_{N \to \infty} \bX^{(N)} (0) / N = \by_0$ for some $\by_0 \in \R_+^k$, the scaled processes $\bX^{(N)} (t) / N$ may be approximated over finite time intervals,  as $N \to \infty$, by the solution $\by (t)$ of the system
\begin{eqnarray}
\frac{d \by}{dt} &=& \sum_{\bl \in \mathcal{L}} \bl \beta_{\bl} ( \by ) 
\label{ODEsystem}
\end{eqnarray}
with $\by(0) = \by_0$ (Theorem~11.2.1 of~\cite{EK05}).
System~\eqref{ODEsystem} is the deterministic epidemic model corresponding to our stochastic model $\bX^{(N)} (t)$.

We assume that the deterministic model~(\ref{ODEsystem}) has two equilibrium points in $\R_+^k$; a stable equilibrium point $\by^*$ (the endemic equilibrium) and an unstable equilibrium point $\by^\circ$  (the disease-free equilibrium).
We further assume that $y^*_i > 0$ for $i=1,2,\ldots,k$, and that $y_i^\circ = 0$ for $i \in \mathcal{I}$, so that $\by^*$ lies in the interior of~$\R_+^k$ while $\by^\circ$ lies on the boundary.
Not all (supercritical) epidemic models fit within this framework, but many standard models do, including all of the illustrative examples that we will present.

Following a successful invasion of infection, the process typically settles into a quasistationary (metastable) endemic phase, before stochastic fluctuations lead to eventual disease extinction.
The time from endemicity to extinction is exponentially distributed~\cite{DP13}, with expected value~$\tau$ satisfying
\begin{eqnarray}
\lim_{N \to \infty} \frac{\ln \tau}{N} &=& U ( \by^\circ ) ,
\label{tau_equation}
\end{eqnarray}
where the function $U ( \by )$  satisfies the Hamilton-Jacobi partial differential equation 
\begin{eqnarray}
H \left( \by , \frac{\partial U}{\partial \by} \right) &=& 0
\label{HJE}
\end{eqnarray}
with $U ( \by^* ) = 0$,
and the Hamiltonian $H ( \by , \btheta )$ is defined to be
\begin{eqnarray}
H ( \by , \btheta ) &=& \sum_{\bl \in \mathcal{L}} \beta_{\bl} ( \by ) \left( {\rm e}^{\bl^T \btheta} - 1 \right ) 
\label{Hamiltonian}
\end{eqnarray}
for $\by \in \R_+^k$, $\btheta \in \R^k$. 

In some special cases (see below), the partial differential equation~(\ref{HJE}) can be solved explicitly for $U ( \by )$.
In general, one must resort to numerical solution, using the method of characteristics (see, for example, section~3.2 of~\cite{E10}).
The characteristic ordinary differential equations (sometimes referred to as Hamilton's equations of motion) are given by
\begin{eqnarray}
\left.
\begin{array}{rcl}
\displaystyle \frac{d\by}{dt} 
&=& \displaystyle \frac{\partial H}{\partial \btheta} 
\;\;=\;\; \displaystyle \sum_{\bl \in \mathcal{L}} \bl \beta_{\bl} ( \by ) {\rm e}^{\bl^T \btheta} ,
\\
\displaystyle  \frac{d\btheta}{dt} 
&=& \displaystyle - \frac{\partial H}{\partial \by} 
\;\;=\;\; \displaystyle  \sum_{\bl \in \mathcal{L}} \frac{d \beta_{\bl}}{d\by} \left( 1 - {\rm e}^{\bl^T \btheta} \right ) .
\end{array}
\right\}
\label{characteristic_equations}
\end{eqnarray}

When $\btheta = {\bf 0}$, equations~(\ref{characteristic_equations}) reduce to the deterministic system~(\ref{ODEsystem}), so that equations~(\ref{characteristic_equations}) have equilibrium points at $( \by , \btheta ) = ( \by^* , {\bf 0} )$ and $( \by^\circ , {\bf 0} )$.
We assume that there also exists a unique equilibrium point of the form $( \by^\circ , \btheta^\circ )$ with $\btheta^\circ \ne {\bf 0}$.
The problem of finding the value of $U ( \by^\circ )$ in equation~(\ref{tau_equation}) then reduces to solving the equations~(\ref{characteristic_equations}) along a characteristic curve~$\Gamma$ connecting the endemic equilibrium point $( \by^* , {\bf 0} )$ at time $t=-\infty$ to the disease-free equilibrium point $( \by^\circ , \btheta^\circ )$ at time $t=+\infty$, and then evaluating
\begin{eqnarray*}
U ( \by^\circ ) &=& \int_{\by^*}^{\by^\circ} \frac{\partial U}{\partial \by} \cdot d\by
\;\;=\;\;
\int_\Gamma \btheta \cdot d\by .
\end{eqnarray*}
That is, we must compute a heteroclinic orbit of the system, the extinction path~$\Gamma$, and then integrate along the extinction path to evaluate~$U ( \by^\circ )$.
The value of $U ( \by^\circ )$ is sometimes referred to as the `action' along the path, and the components of $\btheta = \left( \theta_1 , \theta_2 , \ldots , \theta_k \right)$ as `conjugate variables.'

Before presenting some numerical approaches to solving equations~(\ref{characteristic_equations}), we consider the possibility that equations~(\ref{characteristic_equations}) may be bypassed altogether through direct analytical solution of equation~(\ref{HJE}).

\subsection*{Explicit solution of the Hamilton-Jacobi equation}
For $k=1$ dimensional problems, the Hamilton-Jacobi equation~(\ref{HJE}) reduces to an ordinary differential equation, and it is often possible to rearrange equation~(\ref{HJE}) and integrate to obtain the function~$U (y)$ explicitly.
A well developed theory exists for particular classes of $k=1$ dimensional systems~\cite{AM10,AM17}.
For systems in $k > 1$ dimensions, it is shown in~\cite{C24} that the partial differential equation~(\ref{HJE}) can be solved explicitly for $U ( \by )$ provided certain asymptotic reversibility conditions are satisfied, conditions~(20) and~(21) of~\cite{C24}.
For most $k>1$ dimensional systems of practical interest, however, these conditions are not satisfied, and one must resort to numerical solution of equations~\eqref{characteristic_equations}.

\subsection*{Shooting methods}
Early work on models in $k=2$ dimensions~\cite{DMRH94,HG95,KM08,BM11,AM17} made use of shooting methods (see section~2.4 of~\cite{H07}, chapter~2 of~\cite{K18}, or chapter~16 of~\cite{HW76}).
In this approach, the system~(\ref{characteristic_equations}) is solved (numerically) as an initial value problem, starting from a point close to the endemic equilibrium~$( \by^* , {\bf 0} )$.
The method is considered to succeed if the system evolves to a point sufficiently close to the disease-free equilibrium point~$( \by^\circ , \btheta^\circ )$.
Otherwise, the trial initial point is modified in some appropriate manner~\cite{H07,K18,HW76}, the system~(\ref{characteristic_equations}) solved from this new trial point, and so on until convergence to~$( \by^\circ , \btheta^\circ )$ is achieved.
Shooting methods can work well for $k=2$ dimensional systems (see, for example,~\cite{KM08} and Fig.~1 of \cite{SFBS11}), but generally fail for $k>2$ due to the sensitivity of the extinction path to small perturbations~\cite{FBSS11,SFBS11}.

\subsection*{Finite-time Lyapunov exponents} 
In~\cite{FBSS11,SFBS11}, a method was developed to compute the extinction path by exploiting its sensitivity to perturbations, via the finite-time Lyapunov exponent (FTLE).
The FTLE provides a measure of how sensitively the future behaviour of the system depends upon its current state $( \by , \btheta )$.
It is argued in~\cite{FBSS11,SFBS11} that the extinction path corresponds to a ridge of points having locally maximal FTLE values.
A number of approaches to finding the extinction path via FTLE values have been proposed~\cite{FBSS11,SFBS11,KSS12,BFB16}, but none appear to work well in $k>1$ dimensions.
Consequently, in~\cite{BFB16}, it was proposed to use FTLE values to identify a trajectory reasonably close to the extinction path, which may then be used as the starting point in computing the extinction path more exactly using the `iterative action minimizing method', described below.
Although this approach was successfully implemented in~\cite{BFB16} for particular $k=2$ and $k=3$ dimensional systems, the process involved careful thought and substantial tuning for each specific system.
The FTLE approach does not appear to provide a practical general method for computing extinction paths in $k>1$ dimensions.

\subsection*{Finite difference methods}
In light of the difficulties with shooting and FTLE methods, in~\cite{LS13} a finite difference method was proposed under the name `iterative action minimizing method' (IAMM).
Variants of this approach have since become standard in the literature, e.g.~\cite{BFB16,HS16,NBKFB17,HSS18,KHA22}.

The finite difference approach (see, for example, chapter~2 of~\cite{H07} or chapter~3 of~\cite{K18}) proceeds as follows.
Write $\bz = ( \by , \btheta )$, and consider the system~(\ref{characteristic_equations}) over some finite time interval $\left[ t_s , t_f \right]$ subject to boundary conditions $\bz ( t_s ) = \bz_s$ and $\bz ( t_f ) = \bz_f$, for specified $\bz_s$,~$\bz_f$.
Fix time points $t_1 , t_2 , \ldots, t_{n-1}$ at which the (approximate) solution will be evaluated, with $t_s = t_0 < t_1 < \cdots < t_n = t_f$.
At each point $t_i$, for $i=1 , 2 , \ldots , n-1$, the derivatives $d\bz/dt$ on the left hand side of equations~(\ref{characteristic_equations}) are replaced by an appropriate finite difference approximation~$\bdelta_i$.
A variety of choices for $\bdelta_i$ are available, a few of which are listed in table~1.1 of~\cite{H07}.
For example, if we assume a uniform time step~$h$, with $t_i = t_0 + ih$ for $i=0,1,2,\ldots,n$, then the second order centred finite difference approximation $\bdelta_i$ is given by
\begin{eqnarray}
\bdelta_i
&=&
\frac{\bz_{i+1} - \bz_{i-1}}{2h} .
\label{finite_difference_formula}
\end{eqnarray}

Denoting by $\bz_i = \left( \by_i , \btheta_i \right)$ our approximation to $\bz ( t_i )$ for $i=0,1,2,\ldots,n$, then the ordinary differential equations~(\ref{characteristic_equations}) are thus approximated by the system of algebraic equations
\begin{eqnarray}
\bdelta_i &=& 
\left(
\begin{array}{c}
\sum_{\bl \in \mathcal{L}} \bl \beta_{\bl} ( \by_i ) {\rm e}^{\bl^T \btheta_i} \\ \\
\sum_{\bl \in \mathcal{L}} \left. \frac{d\beta_{\bl}}{d\by} \right|_{\by = \by_i} \left( 1 - {\rm e}^{\bl^T \btheta_i} \right)
\end{array}
\right)
\mbox{ for } i=1,2,\ldots,n-1 .
\label{finite_difference_equations}
\end{eqnarray}
It remains to solve the set of $2k(n-1)$ equations~(\ref{finite_difference_equations}) to find the $2k(n-1)$ unknowns making up the components of $\bz_1 , \bz_2 , \ldots , \bz_{n-1}$.
In~\cite{LS13}, an implementation of Newton's method to iteratively solve equations~(\ref{finite_difference_equations}) is described in detail.

The true extinction path is defined over the interval $t \in [-\infty , \infty]$.
In~\cite{LS13} it is suggested to compute our approximate solution over a finite time interval of the form $\left[ t_s , t_f \right] = [ -T , T]$ for some~$T>0$.
Provided $T$ is sufficiently large, the true extinction path may be expected to stay close to $\left( \by^* , {\bf 0} \right)$ for $-\infty < t < -T$, to move rapidly towards $\left( \by^\circ , \btheta^\circ \right)$ within the transition region $[-T , T]$, and then to remain close to $\left( \by^\circ , \btheta^\circ \right)$ for $T < t < +\infty$.

Having truncated the solution interval to~$[-T , T]$, there are a variety of ways to incorporate the boundary conditions.
Most simply, we may append to equations~(\ref{finite_difference_equations}) the two further equations $\bz_0 = \left( \by^* , {\bf 0} \right)$ and $\bz_n = \left( \by^\circ , \btheta^\circ \right)$, giving a system of $2k (n+1)$~equations as input to the solver.
The boundary conditions are thus treated not as hard constraints, but on an equal footing with equations~(\ref{finite_difference_equations}), so that the solution obtained may be expected to satisfy $\bz_0 \approx \left( \by^* , {\bf 0} \right)$, $\bz_n \approx \left( \by^\circ , \btheta^\circ \right)$.

For a $2k$ dimensional system of first order ordinary differential equations such as~(\ref{characteristic_equations}), over a finite interval~$[ t_s , t_f ]$, one would normally expect to impose a total of $2k$ boundary conditions on the components of $\bz ( t_s )$ and $\bz ( t_f )$.
Indeed, if all $2k$ components of $\bz ( t_s )$ are specified, then the problem may be treated as an initial value problem.
Since we wish to specify both $\bz ( t_s )$ and $\bz ( t_f )$, we have a total of $4k$ boundary conditions, and the system is overdetermined.
This will not necessarily cause any problems, since we expect our system of equations~\eqref{finite_difference_equations} to admit a solution satisfying all boundary conditions. 
We consider other possibilities for the boundary conditions below.

As the extinction path is expected to make a sharp transition near the centre of the domain, it is suggested in~\cite{LS13} to make use of a non-uniform grid~$t_1 , t_2 , \ldots , t_{n-1}$, designed to have higher resolution in the region where the solution is transitioning most rapidly.
The generalisation of the second order centred finite difference formula~(\ref{finite_difference_formula}) to the case of a nonuniform grid is given as formula~(9) of~\cite{LS13}.
We will not pursue this here; the time transformation method described below achieves a similar effect.

The solver needs to be initialised from some initial guess for the values of $\{ \bz ( t_i ) : i=0,1,2,\ldots,n \}$.
One suggestion from~\cite{LS13} is to use as initial guess $\by_i = \by^\circ + ( \by^* - \by^\circ ) / \left( 1 + {\rm e}^{C t_i} \right)$ and $\btheta_i = {\bf 0}$, where $C>0$ is an appropriately chosen constant.
The form of $\by_i$ here is intended to reflect the sharp transition made by the extinction path, with the value of~$C$ adjusting for the sharpness of the transition.
We discuss other options below.

In implementing the finite difference method, the approach of~\cite{LS13} was adapted by~\cite{HSS18} in a number of ways, details of which may be seen in the Matlab code provided as electronic supplemental material to~\cite{HSS18}.
Firstly, eighth order (rather than second order) finite difference formulae are used, with lower order expressions close to the boundaries.
The required higher order finite difference formulae are available from~\cite{F88}.
Secondly, rather than a bespoke implementation of Newton's method, Matlab's {\tt fsolve} function~\cite{fsolve} is used, which simplifies the coding and provides a more flexible solver.
For boundary conditions, the code provided by~\cite{HSS18} appends to equations~(\ref{finite_difference_equations}) the corresponding equations for $i=0$ and $i=n$, but with left hand sides set to zero.
This is appropriate since $\left( \by^* , {\bf 0} \right)$ and $\left( \by^\circ , \btheta^\circ \right)$ are both stationary points of the system.
This choice has the consequence that for any equilibrium point $\bz^\dagger$ of equations~(\ref{characteristic_equations}), the full system of equations supplied to the solver admits the constant solution $\bz_0 = \bz_1 = \cdots = \bz_n = \bz^\dagger$.
In particular, both $\bz_0 = \bz_1 = \cdots = \bz_n = \left( \by^* , {\bf 0} \right)$ and $\bz_0 = \bz_1 = \cdots = \bz_n = \left( \by^\circ , \btheta^\circ \right)$ provide exact solutions to the system of algebraic equations.
Successful performance thus depends upon supplying to the solver an initial guess such that the algorithm will converge towards a solution with $\bz_0 \approx \left( \by^* , {\bf 0} \right)$ and $\bz_n \approx \left( \by^\circ , \btheta^\circ \right)$, and not towards one of these (constant) exact solutions.

Our implementation makes use of ideas from~\cite{HSS18}, but differs in a number of ways, see \nameref{S1_File},~\nameref{S2_File},~\nameref{S3_File} for full details. 
Most significantly, in contrast to~\cite{HSS18}, we `vectorize'~\cite{vectorization} the computation of both sides of equations~\eqref{finite_difference_equations}.
That is, loops are replaced with matrix operations, and functions defined in such a way that they are able to accept multiple input arguments and return the corresponding multiple outputs. 
Because Matlab is optimized for matrix operations, vectorized code can run considerably faster than the corresponding code containing loops~\cite{vectorization}.
Vectorization is particularly worthwhile for our application, because the solver must evaluate the left and right hand sides of equations~(\ref{finite_difference_equations}) repeatedly, and the left hand side of equations~(\ref{finite_difference_equations}) consists of a linear combination of the proposed solution values $\bz_1 , \bz_2 , \ldots , \bz_{n-1}$, so may be computed by a single matrix multiplication.
There is only a small initial set-up cost in creating the (sparse) matrix of coefficients corresponding to the chosen number of grid points~$n+1$ and the specified order of the finite difference approximation.
We found that our vectorized code ran orders of magnitude faster than a non-vectorized version.

\subsection*{Collocation methods}
For numerical solution of boundary value problems, a well-known alternative to finite difference methods, but one that does not seem to have been considered in the current context other than by~\cite{CS24}, is provided by collocation methods.
See, for example, appendix~A.2 of~\cite{K18} or chapter~2 of~\cite{H07} for a general introduction to the approach, which we now briefly outline.

Recall that we aim to solve (approximately) an ordinary differential equation system such as~(\ref{characteristic_equations}) over a finite time interval $[ t_s , t_f ]$ subject to boundary conditions $\bz ( t_s ) = \bz_s$, $\bz ( t_f ) = \bz_f$, where $\bz = \left( \by , \btheta \right)$.
As with the finite difference approach, we first fix time points $t_s = t_0 < t_1 < \cdots < t_n = t_f$ at which to evaluate the solution.
In the collocation approach, the approximate solution $\tilde \bz (t)$ is expressed as
\begin{eqnarray}
\tilde \bz (t) &=& \sum_{j=1}^{2k(n+1)} a_j \bphi_j ( t) ,
\label{collocation_representation}
\end{eqnarray}
where the functions $\bphi_j : \R \to \R^{2k}$ are appropriately chosen basis functions, and~$a_1 , a_2  , \ldots , a_{2k(n+1)}$ are coefficients whose values are to be determined.

Using the representation~(\ref{collocation_representation}), the system~(\ref{characteristic_equations}) may be approximated by the system
\begin{eqnarray}
\sum_{j=1}^{2k(n+1)} a_j \left. \frac{d\bphi_j}{dt} \right|_{t = t_i}
=
\left(
\begin{array}{c}
\sum_{\bl \in \mathcal{L}} \bl \beta_{\bl} ( \by_i ) {\rm e}^{\bl^T \btheta_i} \\ \\
\sum_{\bl \in \mathcal{L}} \left. \frac{d\beta_{\bl}}{d\by} \right|_{\by = \by_i} \left( 1 - {\rm e}^{\bl^T \btheta_i} \right)
\end{array}
\right)
\mbox{ for } i=1,2,\ldots,n-1.
\label{collocation_equations}
\end{eqnarray}

On the right hand side of equation~(\ref{collocation_equations}), recalling that $\bz = \left( \by , \btheta \right)$, each term $\by_i$, $\btheta_i$ may be expressed as a function of $a_1 , a_2 , \ldots , a_{2k(n+1)}$ using the representation~(\ref{collocation_representation}).
To find the values of~$a_1 , a_2 \ldots , a_{2k(n+1)}$, it then remains to solve the algebraic system~(\ref{collocation_equations}) together with the boundary conditions $\bz ( t_s ) = \bz_s$ and $\bz ( t_f ) = \bz_f$, a set of $2k(n+1)$ equations. 

We can deal with the fact that the true extinction path is defined over an infinite time interval exactly as under the finite difference approach, by truncating to the finite interval $[-T , T]$ for suitable~$T>0$.

In practice, we will make use of the Matlab function {\tt bvp5c}, which implements a more sophisticated version of collocation than described above.
For details, see~\cite{KS08}.
An issue that arises here is that when solving the $2k$~dimensional system~(\ref{characteristic_equations}), the {\tt bvp5c} function accepts a total of $2k$ boundary conditions, whereas we would like to impose the $4k$ boundary conditions $\bz_0 =  \bz^*$, $\bz_n = \bz^\circ$, where $\bz^* = \left( \by^* , {\bf 0} \right)$, $\bz^\circ = \left( \by^\circ , \btheta^\circ \right)$.
One way around this would be to impose boundary conditions of the form $\left( \bz_0 - \bz^* \right)_i^2 + \left( \bz_n - \bz^\circ \right)_i^2 = 0$ for $i=1,2,\ldots,2k$.
However, we found that in practice this did not work well.
Instead, we impose the $2k$ boundary conditions $\by_0 = \by^*$, $\btheta_n = \btheta^\circ$.
We then examine the diagnostic values described below (\nameref{convergence_diagnostics}) to check that the end points of the computed extinction path, $\bz_0$,~$\bz_n$, are sufficiently close to the equilibrium points~$\bz^*$,~$\bz^\circ$, respectively.

\subsection*{Convergence diagnostics}\label{convergence_diagnostics}
For each of our algorithms, a number of tuning parameters must be adjusted to obtain satisfactory convergence.
If tuning parameters are poorly chosen, then the code returns an error message and no results.
Even when the code returns results with no error messages, the results may not be reliable, as we shall see in our~\nameref{results} for the network SIS model.

To check convergence, we compute the Euclidean distances between the end points of the computed trajectory and the corresponding equilibrium points of the system~\eqref{characteristic_equations}, as well as the maximal value of the Hamiltonian~\eqref{Hamiltonian} along the computed trajectory. 
That is, we compute the diagnostic values $d^* = || \left( \by_0 , \btheta_0 \right) - \left( \by^* , {\bf 0} \right) ||_2$, $d^\circ = || \left( \by_n , \btheta_n \right) - \left( \by^\circ , \btheta^\circ \right) ||_2$, and $M = \max \left\{ | H \left( \by_i , \btheta_i \right) | : i=0,1,2,\ldots,n \right\}$. 
Note that in our implementations of finite difference methods, the number of grid points~$(n+1)$ remains fixed, whereas in our collocation methods, Matlab's {\tt bvp5c} function automatically adjusts the number of grid points as part of the solution process.
In computing~$M$, to allow for a more direct comparison between methods, we always compute values of $H \left( \by_i , \btheta_i \right)$ at the set of grid points $\{ \left( \by_i , \btheta_i \right) : i=0,1,2,\ldots,n\}$ that we supplied to the finite difference method, and not the adjusted set of grid points generated by {\tt bvp5c}.

Provided the three diagnostic values $d^*$, $d^\circ$, $M$ are all close to zero, then we have a computed path that starts close to the endemic equilibrium point~$\left( \by^* , {\bf 0} \right)$, ends close to the disease-free equilibrium point~$\left( \by^\circ , \btheta^\circ \right)$, and approximately satisfies the Hamilton-Jacobi equation~\eqref{HJE} (with $\btheta = \partial U / \partial \by$) along the entire path.
We can thus have confidence that, however obtained, the computed solution provides a reasonable approximation to the true extinction path.

\subsection*{Parameter continuation}
A critical component of both finite difference and collocation methods is the initial guess, which must be sufficiently close to the true solution for the numerical solver to converge towards the true extinction path.
The most obvious initial guess is a straight line joining the endemic point $( \by^* , {\bf 0} )$ to the disease-free point $( \by^\circ , \btheta^\circ )$, with points equally spaced along the line and a uniform grid in time.
This can work well when the system is close to criticality (that is,  $R_0$ is only slightly above~1), so that the endemic and disease-free points are very close together, but does not work so well for highly supercritical systems ($R_0 \gg 1$).
One way to deal with this is through parameter continuation~\cite{SKR00}.
That is, we solve a sequence of problems, starting with parameter values chosen such that solution is straightforward, and using the solution trajectory for one set of parameter values as initial guess for the problem with slightly modified parameter values.
This process is continued  until the parameter values of interest are attained.
In the epidemic modelling context, the effect of varying parameter values is itself of great interest---in particular, intervention policies can often be modelled through changing the values of model parameters.
Computing solutions across a range of parameter values is thus a very natural approach.

\subsection*{Generating the initial guess from an explicit solution for a related model}
Another way to generate an initial guess for the solution trajectory is to find a model that is closely related to the model of interest, and for which equation~\eqref{HJE} can be explicitly solved for~$U ( \by )$.
The conditions given in~\cite{C24} can help in finding such a model.
The extinction path for the analytical solvable model can then be used to generate an initial guess for the extinction path of the model of interest, as follows.

For the analytically solvable model, knowledge of the solution $U ( \by )$ allows us to write down explicit formulae for the conjugate variables $\boldsymbol{\theta} = \left. \partial U\right/ \partial \boldsymbol{y}$.
Substituting for $\btheta ( \by )$ into equations~(\ref{characteristic_equations}) we obtain a $k$~dimensional system of ordinary differential equations in $\by (t)$.
Provided our analytically solvable model satisfies the conditions~(20) and~(21) of~\cite{C24}, the system thus obtained is precisely the system~(\ref{ODEsystem}) in reversed time.

For many standard epidemic models, including the network SIS model that we will use to illustrate this technique, the system~(\ref{ODEsystem}) is straightforward to solve numerically, because the endemic equilibrium point~$\by^*$ is globally asymptotically stable in the interior of the state space.
Consequently, taking any initial point within the interior of the state space and close to the disease-free equilibrium $\by^\circ = {\bf 0}$, we can numerically solve the resulting initial value problem to obtain a solution trajectory that starts close to~$\by^\circ$ and ends close to the endemic point $\by^*$.
We then reverse in time the trajectory in $\by$ space, append to this the corresponding $\btheta ( \by )$ values from our analytical solution to equation~\eqref{HJE}, and thus obtain a solution to equations~(\ref{characteristic_equations}) for the analytically solvable  model.
This provides the initial guess to supply to the numerical solver used to compute the extinction path for our model of interest.

Numerical solution of equations~(\ref{ODEsystem}) is not quite so straightforward as may at first appear, since we are interested in high dimensional systems, and aim to find a trajectory that starts and ends at equilibrium points.
Consequently, we solve using the Matlab function {\tt ode23s}, designed to deal with stiff differential equation systems, and avoid starting too close to the equilibrium point~$\by^\circ$.
This remains considerably more straightforward than direct numerical solution of equations~(\ref{characteristic_equations}).

\subsection*{Transforming the time interval}
Rather than truncating to the finite time interval $[ - T , T ]$, an alternative, suggested by~\cite{BFB16,CS24}, is to apply a transformation $\hat t = \psi (t)$, where $\psi ( \cdot )$ is a continuously differentiable increasing function mapping $(- \infty , \infty )$ to the finite interval $\left( \hat t_s , \hat t_f \right) = (-T , T)$ for some $T$.
We can then solve (approximately) the transformed version of the system~(\ref{characteristic_equations}) over the interval $\left[ -T , T \right]$ using either a finite difference method or a collocation method.
This has the advantages that (i)~we are now effectively solving over the full (untruncated) time interval $t \in [-\infty , +\infty ]$; and (ii)~a uniform grid of points in transformed time, $\hat t_0 , \hat t_1 , \ldots , \hat t_n$, will correspond to points in untransformed time, $t_0 , t_1 , \ldots , t_n$, that are more widely spaced out as $t \to \pm \infty$, reflecting the nature of the extinction path.

One issue with this approach is that the autonomous system~(\ref{characteristic_equations}) transforms to a non-autonomous system, but this is not a problem in practice, provided computer code is written to allow for explicit time dependence in the derivatives.
A second issue arises at the boundary points, where we have $d\bz/dt = {\bf 0}$, but the Jacobian of the transformation, $d\psi/dt$, will also be zero, and so $d\bz / d\hat t$ is undefined at the boundary.
In the case of the finite difference method, this is straightforward to deal with: instead of appending to equations~(\ref{finite_difference_equations}) the condition that the derivatives be zero at the boundaries, we instead append to the transformed version of equations~(\ref{finite_difference_equations}) the equations $\bz_0 = \left( \by^* , {\bf 0} \right)$ and $\bz_n = \left( \by^\circ , \btheta^\circ \right)$.
It thus becomes unnecessary to evaluate derivatives with respect to~$\hat t$ at the boundaries.
In the case of the collocation method, we simply arrange that the required derivatives evaluate to zero at the boundaries.
This may not be correct, but by examining the diagnostic values $d^*$, $d^\circ$, $M$, we can check that the computed path does, nevertheless, provide a reasonable approximation to the true extinction path.

\section*{Models}\label{models}
We will demonstrate our numerical algorithms using the three epidemic models described below as illustrative examples.

\subsection*{Ross-Macdonald malaria model}
The transmission of malaria between human hosts and mosquito vectors was first modelled mathematically by Ross~\cite{R11}, the model later being developed further by Macdonald~\cite{M57}.
A stochastic version of the Ross-Macdonald model was presented in~\cite{N91}, and various aspects of the model of~\cite{N91} have since been studied by a number of authors~\cite{LZM07,BT17,CS24}.
In particular, the expected extinction time of infection for this model was studied in~\cite{BT17,CS24}.
Although originally developed with malaria in mind, the model can be applied to other vector-borne infections such as dengue fever, yellow fever, and Zika virus disease.

Consider a population consisting of $N$~hosts and $V$~vectors, and set $c = V/N$.
Each individual (whether host or vector) is assumed to be either susceptible to infection, or infected and infectious.
Denote by $X_1 (t)$,~$X_2 (t)$ the numbers of infected hosts and infected vectors, respectively, at time $t \ge 0$, and recall that scaled numbers of individuals $\left( X_1 / N , X_2 / N \right)$ are denoted by $\by = \left( y_1 , y_2 \right)$ in the limit as $N \to \infty$.
Taking transition rates to be of the form~\eqref{density_dependent_rates} with functions $\beta_{\bl} ( \by )$ given in table~\ref{RM_rates}, where $c,\eta, p, q, \sigma , \delta > 0$, we obtain the model studied in~\cite{N91,LZM07,BT17,CS24}.
Here $\eta$~denotes the biting rate of vectors on hosts, $p$~the vector-to-host transmission probability, $q$~the host-to-vector transmission probability, $\sigma^{-1}$~the mean infectious period of hosts, and $\delta^{-1}$ the mean lifetime of infected vectors.

\begin{table}[!ht]
\centering
\caption{
{\bf Transitions and rate functions for the Ross-Macdonald model}}
\begin{tabular}{|l|c|c|}
\hline
Event & Transition vector $\bl$ & Transition rate function $\beta_{\bl} ( \by )$ \\ \hline
Infection of a host & $(1,0)$ & $\eta p \left( 1-y_1 \right) y_2$ \\ \hline
Infection of a vector & $(0,1)$ & $\eta q ( c - y_2 )  y_1$ \\ \hline
Recovery of a host & $(-1 , 0 )$ & $\sigma y_1$ \\ \hline
Death of a vector & $( 0, -1 )$ & $\delta y_2$ \\ \hline
\end{tabular}
\label{RM_rates}
\end{table}

The host-to-host basic reproduction number~$R_0$ for this model, being the expected number of secondary host infections generated by a single infectious host introduced into an otherwise susceptible population, is given by~\cite{LZM07}
\begin{eqnarray}
R_0 &=& \frac{cpq \eta^2}{\sigma\delta} .
\label{R0_Ross-Macdonald}
\end{eqnarray}

For $R_0 > 1$, the endemic and disease-free equilibrium points of the system~(\ref{characteristic_equations}) for this model are~\cite{CS24}, respectively,
\begin{eqnarray*}
\left( \by^* , {\bf 0} \right) 
&=& 
\left( 
\frac{R_0 - 1}{R_0} \left( \frac{cp\eta}{cp\eta+\sigma} \right) ,\ 
\frac{R_0 - 1}{R_0} \left( \frac{cq\eta}{q\eta + \delta} \right) ,\ 
0,\ 0 \right) ,
\end{eqnarray*}
and
\begin{eqnarray*}
\left( \by^\circ , \btheta^\circ \right) 
&=&
\left( 0, \ 0, \ 
\ln \left( \frac{p\eta+\delta}{p\eta+\delta R_0} \right) ,\ 
\ln \left( \frac{cq\eta+\sigma}{cq\eta+\sigma R_0} \right) \right) .
\end{eqnarray*}

\subsection*{Network susceptible-infectious-susceptible (SIS) model}
Models of SIS form have been suggested as appropriate for biological infections that do not induce immunity, such as gonorrhea~\cite{LY76}, as well as for computer viruses spreading through a network~\cite{KW91,KW93,WM04,KHA22}. 
In both cases, it is important to take into account heterogeneous population structure, representing either different types of individuals~\cite{LY76}, or network structure~\cite{C18a}.

Consider a closed population of $N$~individuals divided into $k$~groups, with group~$i$ ($i=1,2,\ldots,k$) consisting of $N_i$~individuals, where $N_1 + N_2 + \cdots + N_k = N$.
Denote by $f_i = N_i / N$ the proportion of the population in group~$i$, and suppose that $f_i > 0$ for $i=1,2,\ldots,k$.
Each individual is assumed to be either susceptible to infection, or infected and infectious.
For $i=1,2,\ldots,k$, denote by $X_i (t)$ the number of infected individuals in group~$i$ at time $t \ge 0$, and recall that the scaled number of individuals $X_i / N$ is denoted by~$y_i$ in the limit as $N \to \infty$.
Taking transition rates to be of the form~\eqref{density_dependent_rates} with functions $\beta_{\bl} ( \by )$ given in  table~\ref{SIS_hetero_rates}, where $\be_i$ denotes the unit vector with $i$th component equal to~1, and assuming that $\beta , \gamma > 0$ and $\lambda_i , \mu_i  > 0$ for $i=1,2,\ldots,k$, we obtain the model studied in~\cite{HS16,C18a,KHA22}.
Here $\beta$ is an overall measure of infectiousness, $\gamma^{-1}$ is the mean infectious period, $\lambda_i$~represents the infectiousness of group~$i$ individuals, and $\mu_i$ represents the susceptibility of group~$i$ individuals.
We assume without loss of generality that $\sum_{i=1}^k \mu_i f_i = \sum_{i=1}^k \lambda_i f_i = 1$.

\begin{table}[!ht]
\centering
\caption{
{\bf Transitions and rate functions for the network SIS model}}
\begin{tabular}{|l|c|c|}
\hline
Event & Transition vector $\bl$ & Transition rate function $\beta_{\bl} ( \by )$ \\ \hline
Infection in group~$i$ & $\be_i$ & $\beta \mu_i ( f_i - y_i ) \sum_{j=1}^k \lambda_j y_j$  \\ \hline
Recovery in group~$i$ & $-\be_i$ & $\gamma y_i$ \\ \hline
\end{tabular}
\label{SIS_hetero_rates}
\end{table}

This model may be interpreted as modelling an infection spreading between individuals connected by an uncorrelated (that is, with no correlations between degrees of neighbouring individuals) directed network, as follows~\cite{C18a,KHA22}.
Set group~$i$ to consist of all individuals having in-degree $d^{\mbox{in}} (i)$ and out-degree $d^{\mbox{out}} (i)$, for all pairs $\left( d^{\mbox{in}} ,\ d^{\mbox{out}} \right)$ that exist in the network.
Denote by $\bar d$ the mean in-degree across the network, noting that this is equal to the mean out-degree.
Denote by $\beta^\prime$ the rate at which infection is transmitted along any edge from an infectious to a susceptible individual.
The rate at which new infections arise in group~$i$ is then
\begin{eqnarray}
\frac{d^{\mbox{in}} (i)}{N \bar d} (N_i - X_i) \sum_{j=1}^k \beta^\prime d^{\mbox{out}} (j) X_j .
\label{network_rate}
\end{eqnarray}
Setting $\beta = \beta^\prime \bar d$, $\mu_i = d^{\mbox{in}} (i) / \bar d$, $\lambda_i = d^{\mbox{out}} (i) / \bar d$, and recalling equation~\eqref{density_dependent_rates}, we see that the expression~\eqref{network_rate} is in agreement with the transition rate function $\beta_{\be_i} ( \by )$ given in table~\ref{SIS_hetero_rates}.
The model thus obtained is known as the `annealed' network approximation~\cite{DGM08}.
Extinction time for the case $\bmu = \blambda$, representing an undirected network, has been previously studied in~\cite{HS16}.

The basic reproduction number~$R_0$ for this model is given by~\cite{C18a}
\begin{eqnarray*}
R_0 &=& \frac{\beta}{\gamma} \sum_{i=1}^k \lambda_i \mu_i f_i .
\end{eqnarray*}

For $R_0 > 1$, defining $D ( \blambda , \bmu )$ to be the unique positive solution of
\begin{eqnarray*}
\frac{\beta}{\gamma} \sum_{i=1}^k \frac{\lambda_i \mu_i f_i}{1 + \mu_i D ( \blambda , \bmu )}
&=& 1 ,
\end{eqnarray*}
then the endemic equilibrium point $\left( \by^* , {\bf 0} \right)$ of the system~(\ref{characteristic_equations}) for this model has components~\cite{C18a}
\begin{eqnarray*}
y_i^* 
&=&
\frac{\mu_i f_i D ( \blambda , \bmu )}{1 + \mu_i D ( \blambda , \bmu )} \mbox{ for } i=1,2,\ldots,k,
\end{eqnarray*}
and the disease-free equilibrium point $\left( \by^\circ , \btheta^\circ \right)$  has $\by^\circ = {\bf 0}$ and~\cite{C18a}
\begin{eqnarray*}
\theta_i^\circ
&=&
- \ln \left( 1 + \lambda_i D ( \bmu , \blambda ) \right) \mbox{ for } i=1,2,\ldots,k .
\end{eqnarray*}

\subsection*{Susceptible-exposed-infectious-removed (SEIR) model}
The susceptible-exposed-infectious-removed model describes the spread of an infection that exhibits a latent period (the `exposed' state) as well as infection-induced immunity (the `removed' state).
Models of SEIR type have been proposed for numerous different infections, including measles~\cite{S84}, mumps~\cite{PMMW23} and Covid~19~\cite{CSBB20}.

Denote by $X_1 (t)$,~$X_2 (t)$,~$X_3 (t)$ the numbers of susceptible, exposed and infectious individuals, respectively, at time $t \ge 0$, denote by~$N$ the typical total population size, and recall that for $i=1,2,3$, the scaled number of individuals $X_i / N$ is denoted by~$y_i$ in the limit as $N \to \infty$.
Note that individuals transition to the `removed' state at the end of their infectious period, but since removed individuals have no influence on further infectious spread, there is no need to keep track of the number of individuals in the `removed' category.
Taking transition rates to be of the form~\eqref{density_dependent_rates} with functions $\beta_{\bl} ( \by )$ given in  table~\ref{SEIR_rates},  where $\beta , \gamma , \nu , \mu > 0$, we obtain the classic SEIR model~\cite{H00}.
Here $\beta$ denotes the infection rate parameter, $\gamma^{-1}$ the mean infectious period, $\nu^{-1}$ the mean latent period, and $\mu^{-1}$ the mean individual lifetime (noting that there is no disease-induced mortality in this model).

\begin{table}[!ht]
\centering
\caption{
{\bf Transitions and rate functions for the SEIR model}}
\begin{tabular}{|l|c|c|}
\hline
Event & Transition vector $\bl$ & Transition rate function $\beta_{\bl} ( \by )$ \\ \hline
Birth of a susceptible individual & $(1, 0 , 0)$ & $\mu$ \\ \hline
Death of a susceptible individual & $(-1, 0 , 0 )$ & $\mu y_1$ \\ \hline
Death of an exposed individual & $(0, -1, 0 )$ & $\mu y_2$ \\ \hline
Death of an infectious individual & $(0, 0, -1 )$ & $\mu y_3$ \\ \hline
Infection & $(-1, 1, 0 )$ & $\beta y_1 y_3$  \\ \hline
End of latent period & $(0, -1, 1 )$ & $\nu y_2$ \\ \hline
Removal & $(0, 0, -1 )$ & $\gamma y_3$ \\ \hline
\end{tabular}
\label{SEIR_rates}
\end{table}

The basic reproduction number~$R_0$ for this model is given by~\cite{H00}
\begin{eqnarray*}
R_0 &=& \frac{\beta \nu}{( \mu + \nu ) ( \mu + \gamma )}.
\end{eqnarray*}

It is straightforward to show that for $R_0 > 1$, the endemic and disease-free equilibrium points of the system~(\ref{characteristic_equations}) for this model are, respectively,
\begin{eqnarray*}
\left( \by^* , {\bf 0} \right)
&=&
\left(
\frac{1}{R_0} , \ 
\frac{\mu}{\mu+\nu} - \frac{\mu (\mu+\gamma)}{\nu \beta} ,\ 
\frac{\nu \mu}{(\mu+\nu) (\mu+\gamma)} - \frac{\mu}{\beta} ,\ 
0,\ 
0,\ 
0
\right) ,
\end{eqnarray*}
and
\begin{eqnarray*}
\left( \by^\circ , \btheta^\circ \right)
&=&
\left(
1,\ 0,\ 0,\ 
0 , \ 
\ln \left( \frac{\beta \mu + (\mu+\nu) (\gamma+\mu))}{\beta (\mu+\nu)} \right) ,\ 
\ln \left( \frac{(\mu+\nu) (\gamma+\mu)}{\beta \nu} \right)
\right) .
\end{eqnarray*}

% Results and Discussion can be combined.
\section*{Results}\label{results}
\subsection*{Ross-Macdonald malaria model}

The Ross-Macdonald model provides a useful initial test case, since it is a low ($k=2$) dimensional model, and for realistic parameter values, the extinction path has a very simple shape.
Baseline parameter values suggested in~\cite{BT17}, with time units of years, are $c=5$, $\eta=73$, $p=0.5$, $q=0.15$, $\sigma^{-1} = 0.014$ and $\delta^{-1} = 0.055$.
Literature references for these values, as appropriate for malaria in a population of human hosts and mosquito vectors, are given in~\cite{BT17}.
With these parameter values, from equation~\eqref{R0_Ross-Macdonald} we have $R_0 \approx 1.54$.
In~\cite{CS24}, the effects of varying model parameter values upon the value of the action integral $U \left( \by^\circ \right)$, and hence upon the expected time to extinction via the relationship~\eqref{tau_equation}, were considered, with each parameter being varied across a biologically plausible range of values; see Fig.~1 of~\cite{CS24}.
Here, our focus is upon the performance of computational algorithms rather than epidemiological interpretation.
We will vary only the biting rate parameter~$\eta$, across a range that goes well beyond the biologically plausible, with other model parameters fixed at their baseline values.

We implemented both the finite difference method and the collocation method, in each case using either a truncated time interval, or a time transformation mapping the real line to a finite interval.
The biting rate was varied from~$\eta = 60$ to~$\eta = 500$, so that the basic reproduction number ranges from~$R_0 = 1.04$ to~$R_0 = 72.2$.
For $R_0$ only slightly above~1, the extinction path is well approximated by a straight line from $\left( \by^* , {\bf 0} \right)$ to $\left( \by^\circ , \btheta^\circ \right)$, so we used this straight line as the initial guess supplied to the solver for $\eta = 60$.
As $\eta$ (and hence $R_0$) grows, this simple initial guess will no longer work directly, so we employed continuation on~$\eta$.
We used $n+1 = 101$ uniformly spaced grid points for our initial guess.

When truncating to a finite time interval~$[-T,T]$, we employed continuation on the truncation parameter~$T$ followed by continuation on~$\eta$.
The transition that the path makes from a neighbourhood of $\left( \by^* , {\bf 0} \right)$ to a neighbourhood of $\left( \by^\circ , \btheta^\circ \right)$ becomes more rapid as~$\eta$ (and hence~$R_0$) grows, and so in order to obtain satisfactory convergence, the value of~$T$ is reduced as $\eta$ increases.

When transforming the time interval, we used the time transformation $\hat t = \psi(t) = \tanh ( C t )$, where, to obtain satisfactory convergence, the scaling factor~$C$ is adjusted as $\eta$ is varied.
For the Ross-Macdonald model, for the parameter values considered, we found that taking $C$ to be proportional to the Euclidean distance between the endemic equilibrium point and the disease-free equilibrium point, $C = C^\prime || \bz^* - \bz^\circ ||_2$, and setting $C^\prime = 2$, worked well.

\begin{figure}[!h]
\includegraphics[scale=1]{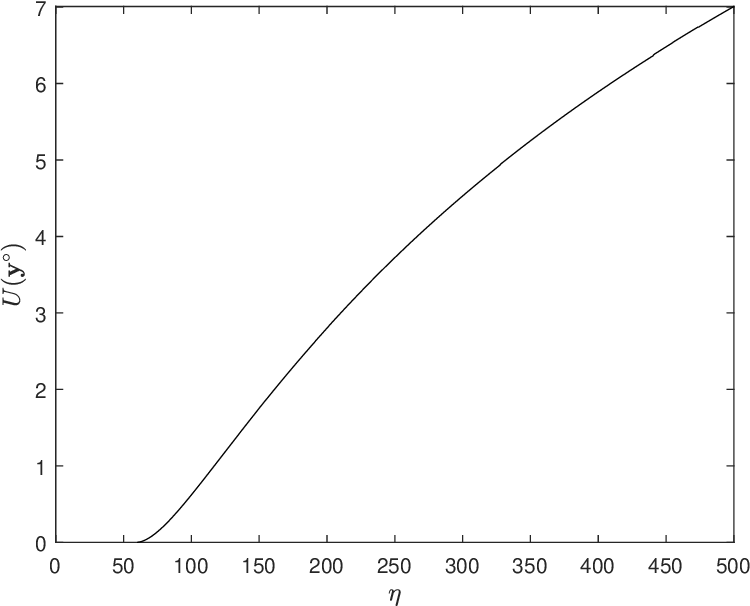}
\caption{{\bf Action integral $U ( \boldsymbol{y}^\circ )$ versus biting rate parameter~$\eta$ for the Ross-Macdonald model.}
Computed using the collocation method with time transformation.
Continuation on~$\eta$ was used, with $\eta$ increasing from~60 to~500, and other model parameters fixed at their baseline values, so that $R_0$ ranges from~1.04 to~72.2.
}
\label{fig_RM_action_versus_eta}
\end{figure}

Fig.~\ref{fig_RM_action_versus_eta} shows the effect of increasing biting rate~$\eta$ upon~$U \left( \by^\circ \right)$, and hence upon mean time to extinction via the relationship~\eqref{tau_equation}.
Note that $U ( \by^\circ )$ is not defined for $R_0 < 1$, and takes the value $U ( \by^\circ ) = 0$ when $R_0 = 1$, and that with other model parameters at their baseline values, $\eta \approx 58.85$ corresponds to $R_0 = 1$.
Fig.~\ref{fig_RM_action_versus_eta} illustrates that interventions that reduce the biting rate (eg bed nets) can be effective in reducing outbreak duration.
For further biological interpretation of plots such as Fig.~\ref{fig_RM_action_versus_eta}, see~\cite{CS24}.

% For figure citations, please use "Fig" instead of "Figure".
% Place figure captions after the first paragraph in which they are cited.

\begin{figure}[!h]
\includegraphics[scale=1]{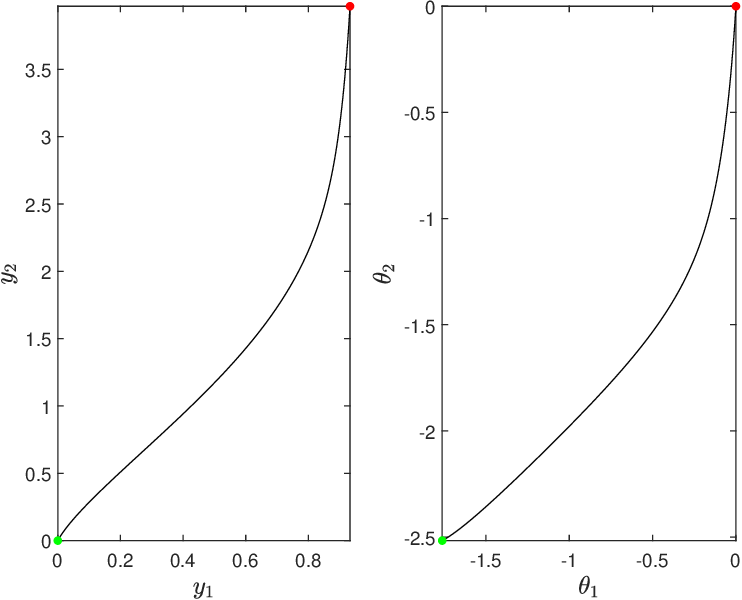}
\caption{{\bf Extinction path for the Ross-Macdonald model.}
Left hand panel shows projection of the 4~dimensional trajectory onto $\left( y_1  , y_2 \right)$ space, right hand panel shows projection onto $\left( \theta_1 , \theta_2 \right)$ space.
Red dot at endemic point $\left( \by^* , {\bf 0} \right)$, green dot at extinction point $\left( \by^\circ , \btheta^\circ \right)$, where $\by^\circ = {\bf 0}$.
Computed using the collocation method with time transformation.
Continuation on~$\eta$ was used to attain the value $\eta=500$, with other model parameters fixed at their baseline values, so that the path shown corresponds to $R_0 = 72.2$.
}
\label{fig_RM_path}
\end{figure}

Fig.\ref{fig_RM_path} shows the computed extinction path for~$\eta = 500$, with other model parameters at their baseline values, so that $R_0 = 72.2$.
For $R_0$ slightly above~1, the extinction path is very close to being a straight line; as $R_0$ increases, the path becomes less linear, but retains a rather simple shape, as illustrated in Fig.~\ref{fig_RM_path}.
Figs.~\ref{fig_RM_action_versus_eta} and~\ref{fig_RM_path} were generated using the collocation method on a transformed time interval; corresponding figures generated by our other three algorithms are essentially indistinguishable.
Note that whereas the finite difference method outputs a solution only at the specified grid points $t_0 , t_1 , \ldots , t_n$, the collocation method, through a representation of the form~(\ref{collocation_representation}), returns a continuous function over the whole solution interval.
The path shown in Fig.~\ref{fig_RM_path} was computed using 45~collocation grid points (selected by the {\tt bvp5c} function), but is plotted based on evaluation of the solution at 1000~equally spaced points, providing a smoother representation of the extinction path.

\begin{figure}[!h]
\includegraphics[scale=1]{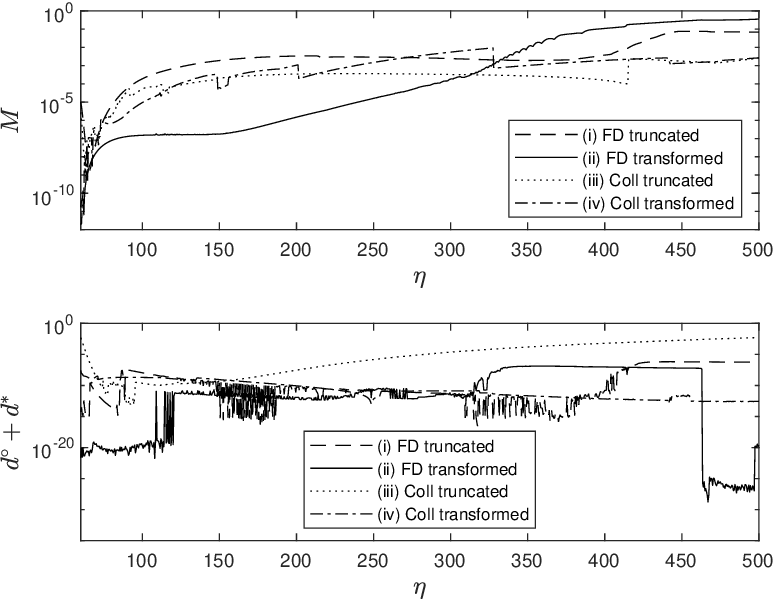}
\caption{{\bf Diagnostic values for the Ross-Macdonald model.}
Upper panel shows maximal absolute value of the Hamiltonian along the computed trajectory, $M = \max \left\{ | H \left( \by_i , \btheta_i \right) | : i=0,1,2,\ldots,n \right\}$.
Lower panel shows sum of Euclidean distances between end points of the computed trajectory and equilibrium points.
Logarithmic scales on vertical axes.
Extinction paths computed using (i)~finite difference method with truncation; (ii)~finite difference method with transformation; (iii)~collocation method with truncation; (iv)~collocation method with transformation.
Continuation on~$\eta$ was used, with $\eta$~increasing from~60 to~500, and other model parameters fixed at their baseline values, so that $R_0$ ranges from~1.04 to~72.2.
}
\label{fig_RM_diagnostics}
\end{figure}

The accuracy of each of the four algorithms is illustrated in Fig.~\ref{fig_RM_diagnostics}, where we see that as~$\eta$ (and hence $R_0$) increases, the error, as measured by $M = \max \left\{ | H \left( \by_i , \btheta_i \right) | : i=0,1,2,\ldots,n \right\}$, tends to increase.
At lower $\eta$ values, the finite difference method with time transformation performs best on this metric.
Collocation methods perform well over a wider range of $\eta$ values than finite difference methods, although poor performance of the finite difference methods only becomes an issue for biologically unrealistic $\eta$ values.
In the lower panel of Fig.~\ref{fig_RM_diagnostics}, for simplicity, we have plotted $d^\circ + d^*$ rather than plotting $d^\circ$ and $d^*$ separately.
We see that in terms of distance between the endpoints of the computed path and the equilibrium points of the system~\eqref{characteristic_equations}, truncation approaches become less accurate than transformation approaches for larger $\eta$ values.

Execution times were as follows: finite difference method with truncation  391~seconds; finite difference method with transformation 701~seconds; collocation method with truncation 13.6~seconds; collocation method with transformation 15.1~seconds.

The diagnostic values illustrated in Fig.~\ref{fig_RM_diagnostics}, and the execution times given above, are heavily influenced by our choices of: number of initial grid points; degree of finite difference formulae; grid of $\eta$~values used in the continuation process; truncation parameter~$T$; transformation function~$\psi(t)$; and transformation scaling parameter~$C$.
The time-consuming part of the process is the experimentation required to find appropriate values for all of the above tuning parameters, rather than the execution time of the final tuned code.
The value of faster execution times is in speeding the process of adjusting tuning parameters.

We have seen that, although accuracy decreases as $R_0$ increases (Fig.~\ref{fig_RM_diagnostics}), for this simple, low-dimensional model we can obtain good results across a very wide range of $R_0$ values using any of our four algorithms.
We now move on to consider a higher dimensional model.
In~\cite{CS24}, an extension of the Ross-Macdonald model that allows for heterogeneity in hosts and vectors was studied.
Fig.~7 of~\cite{CS24} illustrates results for a model with 2~host types and 2~vector types, resulting in a $k=4$ dimensional model.
Rather than pursue this further here, we next consider the network SIS model.

\subsection*{Network SIS model}

The network SIS model provides a useful test case of a higher dimensional model for which, for the parameter values under consideration, the extinction path retains a simple shape.
The extinction path for this model has previously been studied in~\cite{C18a,HS16}.
In Fig.~2 of~\cite{C18a}, results were presented in $k=2$~dimensions for parameter values such that $R_0 = 1.2$.
In Fig.~2(b) of~\cite{HS16}, results were presented for a version of the model in $k=17$~dimensions, with parameter values such that $R_0 = 9$. 
We now consider cases in $k=10$ and $k=30$~dimensions.

We consider an undirected network, and set the degree distribution to be a Poisson distribution of mean~$d$ truncated to the set $\{ 1,2,\ldots,k \}$.
That is, $\bmu = \blambda \propto ( 1,2,\ldots,k )$ with $f_i \propto \frac{d^i}{i!}$ for $i=1,2,\ldots,k$.

In~\cite{C18a}, an analytical solution $U ( \by )$ to the Hamilton-Jacobi equation~(\ref{HJE}) was derived for a directed network in the case~$\blambda = {\bf 1}$, equation~(20) of~\cite{C18a}.
Rather than employing parameter continuation, we use the analytical solution for the case $\blambda = {\bf 1}$ to construct our initial guess for numerical solution in the case $\blambda = \bmu$.
That is, the extinction trajectory was computed independently for each $\beta$ value, using as initial guess the trajectory for the model with the same values of $\beta , \gamma, \boldsymbol{f}$ and $\bmu$, but with $\blambda = {\bf 1}$.

Consider first the case~$k=10$.
We fixed $d=4$, $\gamma = 1$, and solved for a range of values of the overall infection rate parameter~$\beta$, from $\beta=2$ to $\beta =12$, corresponding to $R_0$ ranging from~2.44 to~14.65.
We used $n+1 = 101$ uniformly spaced grid points for our initial guess.
In the same way as for the Ross-Macdonald model, when solving over a truncated time interval $[-T,T]$, the value of $T$ is reduced as $\beta$ is increased.
When solving over a transformed time interval, we again used the transformation $\hat t = \psi(t) = \tanh (Ct)$, but now we set $C = C^\prime || \bz^* - \bz^\circ ||_2^2$ with $C^\prime = 0.1$.
That is, rather than taking the scaling factor~$C$ to be proportional to the distance between the endemic equilibrium point and the disease-free equilibrium point, we found it more effective here to take $C$ to be proportional to the square of the distance.

\begin{figure}[!h]
\includegraphics[scale=1]{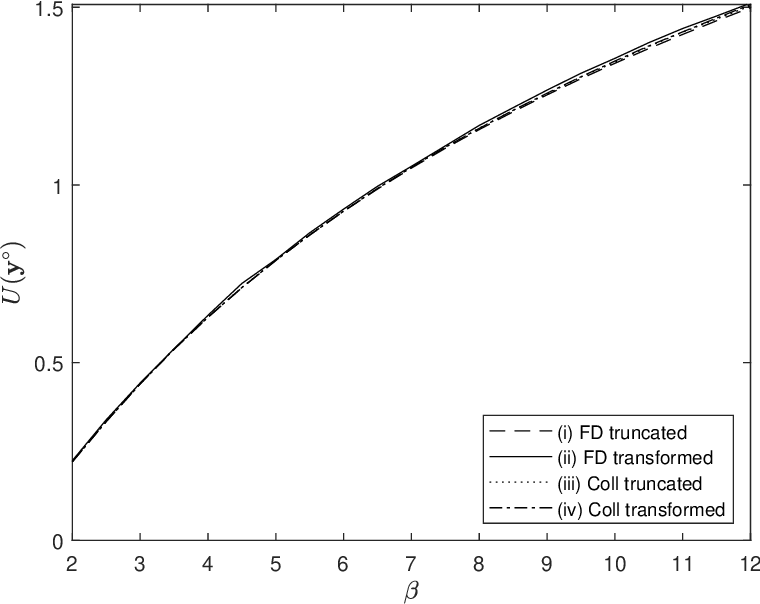}
\caption{{\bf Action integral $U ( \boldsymbol{y}^\circ )$ versus overall infection rate parameter $\beta$ for the network SIS model.}
Computed using (i)~finite difference method with truncation; (ii)~finite difference method with transformation; (iii)~collocation method with truncation; (iv)~collocation method with transformation.
Initial guesses for the extinction path generated from the analytically solvable model with $\blambda = {\bf 1}$.
As $\beta$ is varied from~2 to~12, with other model parameters fixed at values set out in main text, $R_0$ ranges from~2.44 to~14.65.
}
\label{fig_SIS_action_versus_beta}
\end{figure}

Fig.~\ref{fig_SIS_action_versus_beta} shows the effect of increasing the overall infection rate parameter~$\beta$ (and hence $R_0$) upon~$U ( \by^\circ )$.
We see that results obtained by our four algorithms are very similar, although there are small discrepancies, and these discrepancies increase as $\beta$~increases.
We see that the value of $U ( \by^\circ )$, and hence the expected extinction time, is an increasing function of the overall infection rate parameter~$\beta$, as one would expect.
Thus interventions that reduce the overall infection rate parameter can be an effective way to reduce outbreak duration, and Fig.~\ref{fig_SIS_action_versus_beta} allows us to quantify the effect, via the relationship~\eqref{tau_equation}.
	
\begin{figure}[!h]
\includegraphics[scale=1]{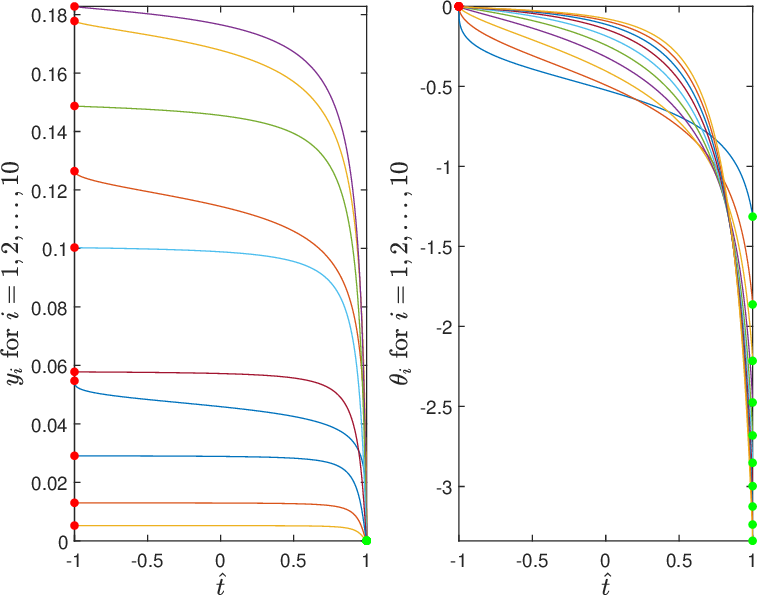}
\caption{{\bf Extinction path for the network SIS model with $k=10$.}
Left hand panel shows state variables $\left( y_1  , y_2  , \ldots , y_{10} \right)$ plotted against transformed time~$\hat t$, right hand panel shows conjugate  variables $\left( \theta_1 , \theta_2 , \ldots , \theta_{10} \right)$ plotted against transformed time~$\hat t$.
Red dots show the components of the endemic point $\left( \by^* , {\bf 0} \right)$, green dots show the components of the extinction point $\left( \by^\circ , \btheta^\circ \right)$, where $\by^\circ = {\bf 0}$.
Computed using the collocation method with time transformation.
Initial guess for the extinction path generated from the analytically solvable model with $\blambda = {\bf 1}$.
Overall infection rate parameter $\beta=12$, with other model parameters values as set out in main text, so that that the path shown corresponds to $R_0 = 14.65$.
}
\label{fig_SIS_path}
\end{figure}

Fig.~\ref{fig_SIS_path} shows the computed extinction path for $\beta = 12$ (so $R_0 = 14.65$).
The solution shown here was computed using the collocation method on a transformed time interval.
The extinction path is a path in $2k = 20$ dimensional space, so we have plotted each of the variables $y_i$, $\theta_i$, $i=1,2,\ldots,10$, against transformed time $\hat t$.

\begin{figure}[!h]
\includegraphics[scale=1]{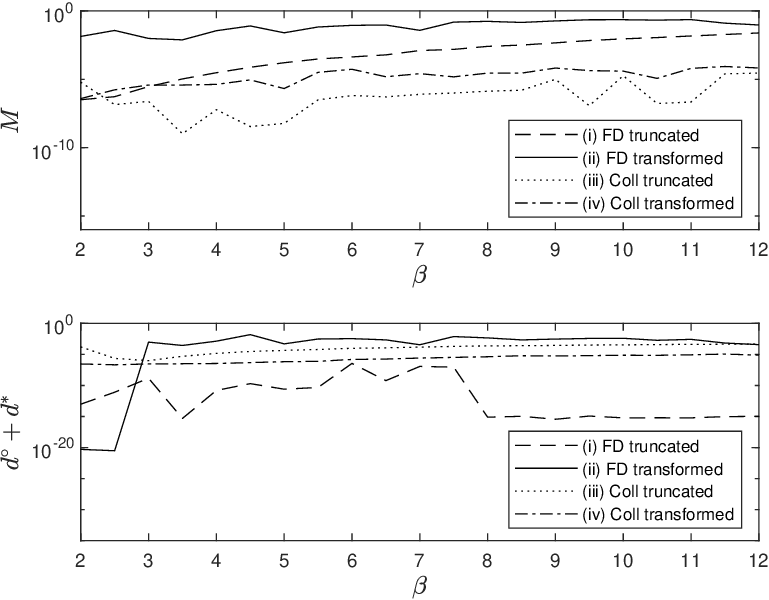}
\caption{{\bf Diagnostic values for the SIS model in a heterogeneous population.}
Upper panel shows maximal absolute value of the Hamiltonian along the computed trajectory, $M = \max \left\{ | H \left( \by_i , \btheta_i \right) | : i=0,1,2,\ldots,n \right\}$.
Lower panel shows sum of Euclidean distances between end points of the computed trajectory and equilibrium points.
Logarithmic scales on vertical axes.
Extinction paths computed using (i)~finite difference method on a truncated time interval; (ii)~finite difference method with a time transformation; (iii)~collocation method on a truncated time interval; (iv)~collocation method with a time transformation.
As $\beta$ is varied from~2 to~12, with other model parameters fixed at values set out in main text, $R_0$ ranges from~2.44 to~14.6.
Initial guesses for the extinction path generated from the analytically solvable model with~$\blambda = {\bf 1}$.
}
\label{fig_SIS_diagnostics}
\end{figure}

Fig.~\ref{fig_SIS_diagnostics} shows diagnostic values for our four algorithms.
In terms of the maximal value of the Hamiltonian along the computed trajectory, $M$, we see that here collocation methods perform better than  finite difference methods, and that truncation of the time interval performs better than time transformation.
In terms of distance between the endpoints of the computed path and the equilibrium points of the system~\eqref{characteristic_equations}, the finite difference method with truncation of the time interval performs best over most of the range of $\beta$ values, and the finite difference method with time transformation performs worst, with collocation methods giving intermediate results.
There is some tendency to loss of accuracy (increasing $M$ values) as $\beta$ (and hence $R_0$) increases, though not so clearly as seen for the Ross-Macdonald model in Fig.~\ref{fig_RM_diagnostics}.

Execution times for the model with $k=10$ were as follows: finite difference method with truncation 256~seconds; finite difference method with transformation 654~seconds; collocation method with truncation 14.1~seconds; collocation method with transformation 541~seconds.

In summary, all four algorithms gave reasonably satisfactory results for $k=10$, with collocation methods found to be more accurate (smaller $M$ values) than finite difference methods. 

\begin{figure}[!h]
\includegraphics[scale=1]{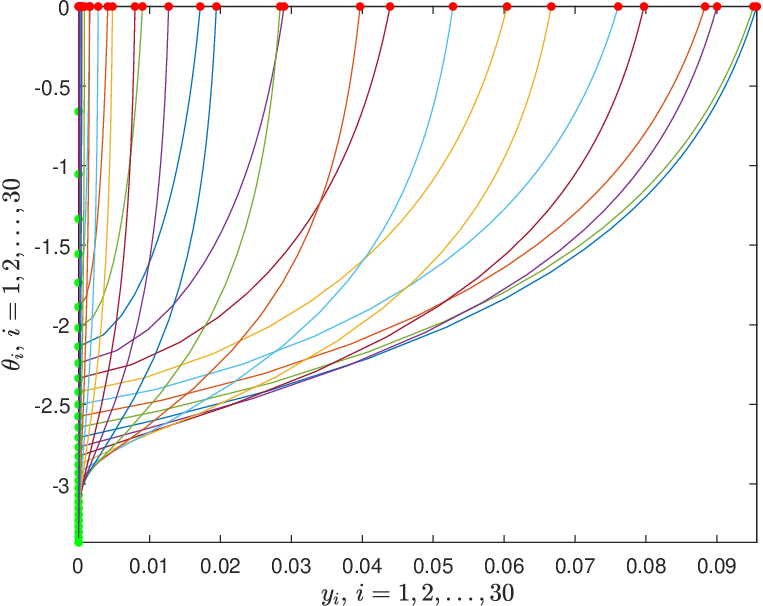}
\caption{{\bf Extinction path for the network SIS model with $k=30$.}
Shown are the 2~dimensional projections $\left( y_i , \theta_i \right)$, for $i=1,2,\ldots,30$, of the 60~dimensional extinction path.
Red dots show the components of the endemic point $\left( \by^* , {\bf 0} \right)$, green dots show the components of the extinction point $\left( \by^\circ , \btheta^\circ \right)$, where $\by^\circ = {\bf 0}$.
Computed using the collocation method with time transformation.
Initial guess for the extinction path generated from the analytically solvable model with $\blambda = {\bf 1}$.
Path shown corresponds to $R_0 = 16$; individual model parameter values set out in main text. 
}
\label{fig_SIS_30_path}
\end{figure}

We next consider the case~$k=30$.
As we increase the dimensionality of the model, it becomes harder to obtain satisfactory convergence, requiring careful adjustment of the various tuning parameters.
Consequently, we present only results obtained using the collocation method with time transformation.
We solved for one set of parameter values, with $d=15$, $\beta = 15$ and $\gamma = 1$, so that $R_0 = 16$.
Fig.~\ref{fig_SIS_30_path} shows the computed extinction path.
As an alternative to the format of Fig.~\ref{fig_SIS_path}, the $2k = 60$ dimensional extinction path is depicted in Fig.~\ref{fig_SIS_30_path} by plotting the 2~dimensional projections $( y_i , \theta_i )$ for $i=1,2,\ldots,30$.
The value of the action integral was computed to be $U ( \by^\circ ) = 1.754$, with diagnostic values $M = 6.76 \times 10^{-6}$, $d^* =  4.61 \times 10^{-5}$, $d^\circ = 2.7 \times 10^{-7}$, suggesting that satisfactory convergence has been achieved.
Execution time was 242~seconds.

The $R_0$ values that we have considered are larger than would be encountered in practice for infections of SIS type.
In general, it becomes harder to obtain satisfactory convergence as $R_0$ increases, so these large $R_0$ values were chosen to test the limits of our algorithms.
We have shown that our four algorithms can all produce satisfactory results for a high dimensional model, even for large $R_0$ values, provided the shape of the extinction path is simple.
We have also demonstrated that generating an initial solution guess from an associated, analytically solvable, model can provide a practical alternative to parameter continuation.
We now move on to consider a model for which the extinction path has a more complicated shape (the SEIR model), making numerical computation considerably more challenging. 

\subsection*{SEIR model}
The extinction path of the SEIR model has previously been studied in~\cite{BFB16}.
In Fig.~14 of~\cite{BFB16}, the extinction path is shown for the case $\mu=0.2$, $\nu = 35$, $\gamma = 100$, $\beta=105$, so that $R_0 = 1.042$.
We fix~$\mu$, $\nu$, $\gamma$ at these same values, but consider a range of values for the infection rate parameter~$\beta$.

When using a time transformation, we found that for this model the transformation $\psi (t) = \tanh ( C t )$ did not produce satisfactory results, so instead we transformed using the `error function',
\begin{eqnarray*}
\psi (t) &=& \mbox{erf} ( Ct )
\;\;=\;\; \frac{2}{\sqrt{\pi}} \int_0^{Ct} \exp ( - s^2 ) \, ds ,
\end{eqnarray*}
for some~$C>0$.
The best results that we were able to achieve were with $C = C^\prime || \bz^* - \bz^\circ ||_2^{0.4}$, where $C^\prime = 0.34$.

For the smallest $\beta$ value considered, we used a straight line from $\left( \by^* , {\bf 0} \right)$ to $\left( \by^\circ , \btheta^\circ \right)$ as the initial guess supplied to the solver.
For all methods, we used parameter continuation on~$\beta$; when truncating the time interval, we first employed continuation on the truncation parameter~$T$.

For collocation methods, we used $n+1 = 301$ uniformly spaced grid points for our initial guess, and computed extinction paths for $\beta$ values from $\beta = 102$ up to $\beta = 160$ in steps of~$0.05$.
For finite difference methods, we used $n+1 = 1201$ uniformly spaced grid points, and $\beta$ values from $\beta = 101$ up to $\beta = 160$ in steps of~$0.01$.
These differences arose because the finite difference methods had difficulty computing the extinction path for $\beta = 102$ directly, so we started the continuation process a little closer to the criticality threshold ($\beta =101$ corresponds to $R_0 = 1.0023$); similarly, a spacing of $0.05$ between consecutive $\beta$ values proved too great for the finite difference methods; and finite difference methods required more than~301 grid points to produce satisfactory results.
Our largest $\beta$ value, $\beta = 160$, corresponds to $R_0 = 1.588$.

\begin{figure}[!h]
\includegraphics[scale=1]{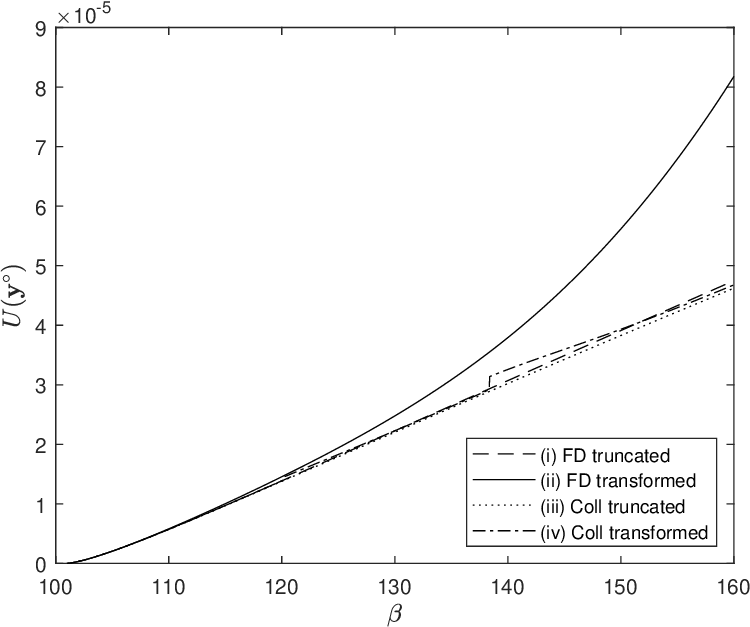}
\caption{{\bf Action integral~$U ( \boldsymbol{y}^\circ )$ versus infection rate parameter~$\beta$ for the SEIR model.}
Computed using (i)~finite difference method with truncation; (ii)~finite difference method with transformation; (iii)~collocation method with truncation; (iv)~collocation method with transformation.
Continuation on~$\beta$ was used, with other model parameters fixed at values set out in main text.
As $\beta$ is varied from~101 to~160, $R_0$ ranges from~1.0023 to~1.588.
}

\label{fig_SEIR_action_versus_beta}
\end{figure}

Fig.~\ref{fig_SEIR_action_versus_beta} shows values of the action integral $U \left( \by^\circ \right)$ computed by each of our four algorithms across a range of $\beta$ values. 
In contrast to previous examples, our four algorithms produced noticeably different results.
All four algorithms are in good agreement up to around $\beta = 120$, but as $\beta$ increases beyond this, results from the finite difference method with time transformation start to diverge strongly from results obtained by other approaches.
Values of $U \left( \by^\circ \right)$ computed using the collocation method with time transformation are in good agreement with other algorithms up to $\beta = 138.35$, but then make a jump upwards at $\beta = 138.4$, before moving back towards results obtained from other algorithms as $\beta$ continues to increase.

\begin{figure}[!h]
\includegraphics[scale=1]{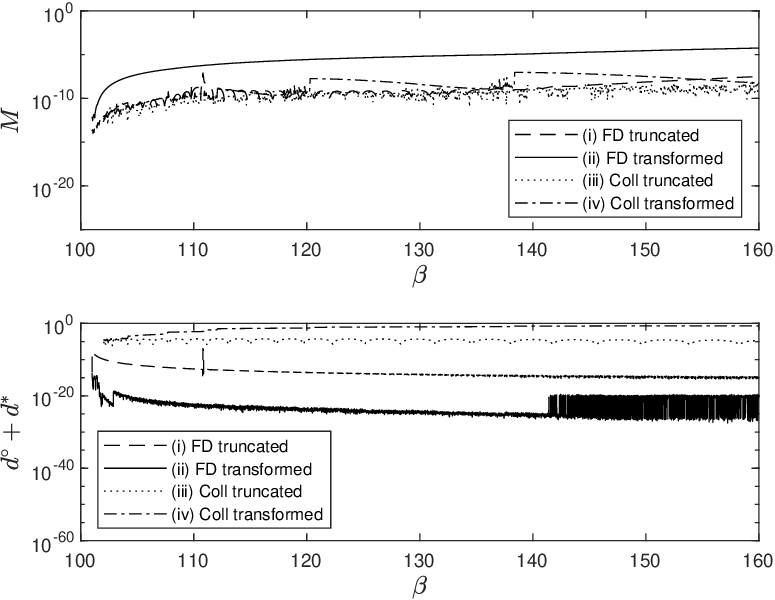}
\caption{{\bf Diagnostic values for the SEIR model.}
Upper panel shows maximal absolute value of the Hamiltonian along the computed trajectory, $M = \max \left\{ | H \left( \by_i , \btheta_i \right) | : i=0,1,2,\ldots,n \right\}$.
Lower panel shows sum of Euclidean distances between end points of the computed trajectory and equilibrium points.
Logarithmic scales on vertical axes.
Extinction paths computed using (i)~finite difference method with truncation; (ii)~finite difference method with transformation; (iii)~collocation method with truncation; (iv)~collocation method with transformation.
Continuation on $\beta$ was used.
As $\beta$ is varied from~101 to~160, with other model parameter values fixed at values set out in main text, $R_0$ ranges from~1.0023 to~1.588.
}
\label{fig_SEIR_diagnostics}
\end{figure}

Fig.~\ref{fig_SEIR_diagnostics}. shows diagnostic values for our four algorithms.
As $\beta$ (and hence $R_0$) increases, the maximal value of the Hamiltonian along the computed trajectory,~$M$, become considerably larger for the finite difference method with time transformation than for the other three algorithms.
It seems clear that results from the finite difference method with time transformation are not reliable here.
We also see an upwards jump in the value of~$M$ for the collocation method with time transformation at $\beta = 138.4$, followed by a gradual decrease as $\beta$ increases further, consistent with the behaviour seen in Fig.~\ref{fig_SEIR_action_versus_beta}, 
This behaviour demonstrates that even when one of our algorithms produces inaccurate results for certain parameter values, computed trajectories may nevertheless converge back towards the true extinction path as the parameter continuation process progresses.

The two methods with time truncation are in reasonably good agreement with one another across the full range of $\beta$ values (Fig.~\ref{fig_SEIR_action_versus_beta}), and have quite similar values of~$M$ for $\beta$ values up to around $\beta = 140$ (Fig.~\ref{fig_SEIR_diagnostics}).
At the largest $\beta$ values, close to $\beta = 160$, values of $M$  start to grow larger for the finite difference method with time truncation than for the collocation method with time truncation (Fig.~\ref{fig_SEIR_diagnostics}), suggesting that overall, the collocation method with time truncation is to be preferred here.

Examining the lower panel of Fig.~\ref{fig_SEIR_diagnostics}, we see that collocation methods do not perform as well as finite difference methods in terms of the diagnostic $d^\circ + d^*$.
That is, using collocation methods, the end points of the computed extinction path are not as close to the equilibrium points of the system~\eqref{characteristic_equations} as we might like. 
However, even on this metric, the performance of the collocation method with time truncation seems acceptable.

\begin{figure}[!h]
\includegraphics[scale=1]{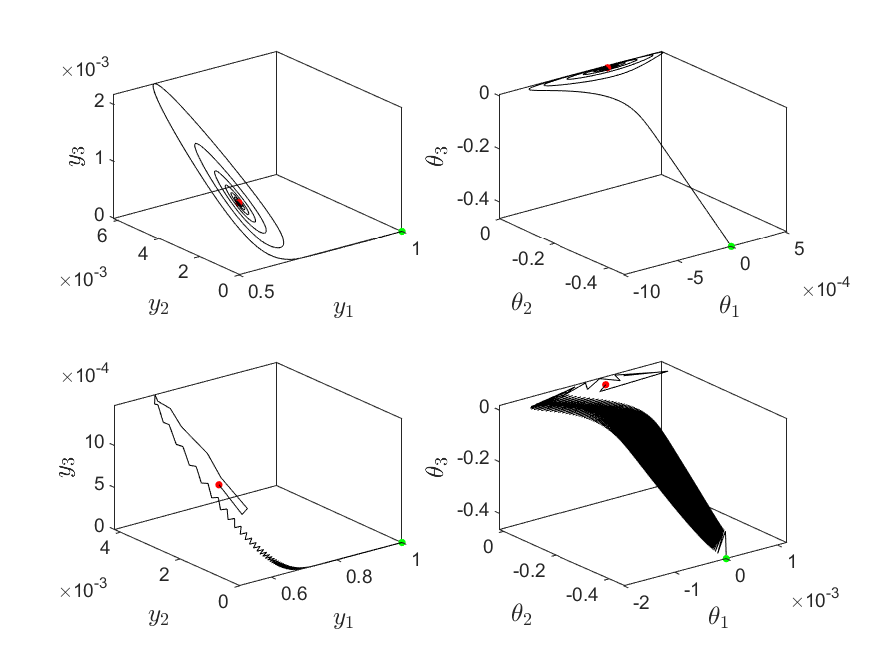}
\caption{{\bf Extinction path for the SEIR model.}
Upper panels show path computed using the collocation method with time truncation; lower panels show path computed using the finite difference method with time transformation.
Left hand panels show projection into $\left( y_1  , y_2  , y_3 \right)$ space; right hand panels show projection into $\left( \theta_1 , \theta_2 , \theta_3 \right)$ space.
Red dots show the endemic point~$\left( \by^* , {\bf 0} \right)$, green dots show the extinction point~$\left( \by^\circ , \btheta^\circ \right)$.
Continuation on~$\beta$ was used to attain the value $\beta = 160$, with other model parameter values fixed at values set out in main text, so that the path shown corresponds to $R_0 = 1.588$.
}
\label{fig_SEIR_path}
\end{figure}

Fig.~\ref{fig_SEIR_path} shows the extinction path computed by the collocation method with time truncation (upper panels) and by the finite difference method with time transformation (lower panels), each with~$\beta = 160$, so that $R_0 = 1,588$.
We see that the path computed by the collocation method with time truncation is a smooth curve, whereas the path computed by the finite difference method with time transformation exhibits jagged zig-zags.
This appears to confirm the evidence of Fig.~\ref{fig_SEIR_action_versus_beta} and Fig.~\ref{fig_SEIR_diagnostics}, that the finite difference method with time transformation is not providing reliable results here.
The collocation method with time truncation, on the other hand, appears to provide reliable results across the full range of $\beta$ values considered.

Execution times were as follows: finite difference method with truncation 68682~seconds; finite difference method with transformation 68030~seconds; collocation method with truncation 1213~seconds; collocation method with transformation 1888~seconds.
We observe that the method producing the most (apparently) reliable results (collocation with truncation) also has the fastest execution time.

% For figure citations, please use "Fig" instead of "Figure".
% Place figure captions after the first paragraph in which they are cited.

% Place tables after the first paragraph in which they are cited.

%PLOS does not support heading levels beyond the 3rd (no 4th level headings).

\section*{Discussion}\label{discussion}
Computing the extinction path for epidemic models is known to be a difficult problem, due to the sensitivity of the path to small perturbations~\cite{FBSS11,SFBS11,BFB16}.
The methods and code that we have presented represent a step towards making WKB methodology accessible, although much scope remains for further progress.

Convergence of the numerical algorithms becomes more difficult to achieve as the basic reproduction number $R_0$ increases far above~1; as dimensionality increases; and when the shape of the extinction path is complicated.
The general tendency to loss of accuracy as $R_0$ increases is illustrated through the diagnostic~$M$ values shown in the upper panels of Fig.s~\ref{fig_RM_diagnostics},~\ref{fig_SIS_diagnostics} and~\ref{fig_SEIR_diagnostics}, while the difficulties in computing extinction paths of more complicated shape are demonstrated in Fig.~\ref{fig_SEIR_action_versus_beta} and the lower panels of Fig.~\ref{fig_SEIR_path},

For each of the three illustrative examples that we have considered, we have been able to present results that go well beyond previously available results for these models.
For the Ross-Macdonald model, in~\cite{CS24} $R_0$ values up to $R_0  = 1.92$ were considered; we have presented results up to $R_0 = 72.2$.
For the network SIS model, in~\cite{HS16} results were presented for an example in $k=17$ dimensions with $R_0 = 9$; we have presented results in  $k=30$ dimensions with $R_0 = 16$.
For the SEIR model, in~\cite{BFB16} results were presented for~$R_0 = 1.042$; we have presented results up to $R_0 = 1.588$.
In addition, for each model, we have computed extinction paths and corresponding values of the action integral~$U ( \by^\circ )$ across whole ranges of parameter values, which is valuable when studying the effects of intervention policies whose effects may be modelled via changes in model parameter values.
Previous authors have often presented results from WKB methodology for only one set of parameter values at a time, e.g.~\cite{BFB16,NBKFB17}.

We have seen that all four of the algorithms considered can work well, even in high dimensions and for large $R_0$ values, provided the extinction path has a simple shape.
For both the Ross-Macdonald model and the network SIS model we were able to compute extinction paths even for $R_0$ values much larger than those of practical relevance.
For the Ross-Macdonald model applied to malaria, baseline parameter values given in~\cite{BT17} correspond to an estimated $R_0$ value of~$1.54$.
For gonorrhea, which may be modelled using the network SIS model, $R_0$ has been estimated to be around~$1.18$--$2$~\cite{BNPS94}.

When the shape of the extinction path is more complicated, the situation is less satisfactory.
For the SEIR model, in Fig.~\ref{fig_SEIR_action_versus_beta} we presented results only up to $R_0 = 1.588$, and even within that range, only two of our four algorithms gave apparently satisfactory results.
Infections to which the SEIR model may be applied include, for example, measles, for which $R_0$ has been estimated to be in the range 12--18~\cite{GMLH17}, and Covid~19, for which $R_0$ has been estimated to be in the range 2.9--9.5~\cite{LR22}.

The shape of the extinction path for the SEIR model is typical of many epidemic models.
Specifically, the path exhibits oscillatory behaviour around $\by^*$ (which lies in the interior of the state space) but not around $\by^\circ$ (on the boundary of the state space).
Other models displaying this pattern include the classic SIR model~\cite{KM08,SFBS11} and the Ebola model of~\cite{NBKFB17}.
Further progress towards reliably computing extinction paths of this shape for larger $R_0$ values would therefore be very valuable.

One natural direction for further investigation is the use of different time transformation functions $\psi (t)$, and in particular, asymmetric functions.
The transformations that we have used, $\psi (t) = \tanh ( Ct )$ and $\psi (t) = \mbox{erf} (Ct)$, are both symmetrical about $t=0$, whereas the shape of the SEIR extinction path (Fig.~\ref{fig_SEIR_path}) is  highly asymmetrical.
Choice of the transformation scaling parameter~$C$ is another aspect that could be investigated further.
We have taken $C$ to be a function of the distance between the equilibrium points $\by^*$ and $\by^\circ$, and found that this worked well.
For other models, it may prove more effective to take $C$ to be a function of $R_0$, or indeed any function of the model parameters.

The use of an analytical solution from a closely related model to generate the initial guess supplied to the solver is an idea that could be further explored.
Our application of this idea to the network SIS model is, so far as we are aware, the first use of this approach.
The recent results of~\cite{C24} open up the possibility of further progress here.

Finally, when employing parameter continuation, we took values for the continuation parameter to be uniformly spaced.
It may be possible to improve convergence through a more sophisticated choice of values of the continuation parameter.
Ideas such as those presented in~\cite{DF89,FD91} suggest that this may be a productive direction for further work.

%\section*{Conclusion}

\section*{Supporting information}

% Include only the SI item label in the paragraph heading. Use the \nameref{label} command to cite SI items in the text.
%\paragraph*{S1 Fig.}
%\label{S1_Fig}
%{\bf Bold the title sentence.} Add descriptive text after the title of the item (optional).

\paragraph*{S1 File.}
\label{S1_File}
{\bf Matlab code for the Ross-Macdonald model.}

\paragraph*{S2 File.}
\label{S2_File}
{\bf Matlab code for the network SIS model.}

\paragraph*{S3 File.}
\label{S3_File}
{\bf Matlab code for the SEIR model.}

\paragraph*{S1 Appendix.}
\label{S1_Appendix}
{\bf Sketch justification of the relationship~\eqref{tau_equation}.}

Recall that we consider a sequence of Markov processes $\{ \bX^{(N)} (t) : t \ge 0 \}$ on $\Z_+^k$, with state space~$S = C \cup \bar C \subseteq \Z_+^k$, such that the process will hit $\bar C$ (the disease-free states) within finite time with probability~1.
Transition rates are of the form~\eqref{density_dependent_rates}, and the scaled processes may be approximated over finite time intervals, for large~$N$, by the solution $\by (t)$ of the system~\eqref{ODEsystem}.
The system~\eqref{ODEsystem} has two equilibrium points in $\R_+^k$; a stable equilibrium point $\by^*$ (the endemic equilibrium)  and an unstable equilibrium point $\by^\circ$ (the disease-free equilibrium).

We assume that for each $N$, there exists a unique quasistationary distribution $\bu^{(N)} = \{ u_{\bx}^{(N)} : \bx \in C \}$ such that for every $\bx , \bx_0 \in C$, 
\begin{eqnarray*}
u_{\bx}^{(N)} &=& 
\lim_{t \to \infty} 
\Pr \left( \bX^{(N)} (t) = \bx \mid \bX^{(N)} (0) = \bx_0 ,\ \bX^{(N)} (t) \in C \right) .
\end{eqnarray*}
This quasistationary distribution represents the metastable behaviour of the process during the endemic phase.

It is known~\cite{DP13} that the time to extinction from quasistationarity is exponentially distributed, and the quasistationary distribution~$\bu^{(N)}$ and expected time to extinction from quasistationarity,~$\tau^{(N)}$, satisfy the quasistationary Kolmogorov forward equation
\begin{eqnarray}
\sum_{\bl \in \mathcal{L}} 
\left( u_{\bx - \bl}^{(N)} \beta_{\bl} \left( \frac{\bx-\bl}{N} \right) - u_{\bx}^{(N)} \beta_{\bl} \left( \frac{\bx}{N} \right) \right)
&=&
- ( \tau^{(N)} N )^{-1} u_{\bx}^{(N)}
\mbox{ for } \bx \in C ,
\label{balance_equation}
\end{eqnarray}
with
\begin{eqnarray}
\tau^{(N)} 
&=&
\left( N \sum_{\bx \in \bar C} \sum_{\bl \in \mathcal{L}} u_{\bx-\bl}^{(N)} \beta_{\bl} \left( \frac{\bx-\bl}{N} \right) \right)^{-1} .
\label{tau_formula}
\end{eqnarray}

From equation~(\ref{balance_equation}), the expected time from endemicity (quasistationarity) to extinction may be found as an eigenvalue of the transition rate matrix of the process.
However, direct evaluation of this eigenvalue is generally not feasible in practice, due to the size of the transition rate matrix, hence we proceed as follows.
Adopting the WKB (Wentzel, Kramers, Brillouin) {\em ansatz}~\cite{AM10,AM17}, we seek a solution of the form
\begin{eqnarray}
u_{\bx}^{(N)} 
&=&
\exp \left( - N U ( \bx / N ) + o(N) \right) \mbox{ as } N \to \infty
\label{WKB}
\end{eqnarray}
for some function $U : \R_+^k \to \R$ that does not depend upon $N$.

Assuming that $\tau^{(N)}$ is sufficiently large for the right hand side of equation~(\ref{balance_equation}) to be neglected, substituting from~(\ref{WKB}) into equation~(\ref{balance_equation}), and collecting leading order terms, then with $H ( \by , \theta )$ defined by equation~\eqref{Hamiltonian}, we find that $U ( \by )$ satisfies the Hamilton-Jacobi equation~\eqref{HJE} with $U ( \by^* ) = 0$.
Assuming that the sum on the right hand side of~(\ref{tau_formula}) is dominated by terms corresponding to states~$\bx-\bl$ in a small neighbourhood of $N \by^\circ$, the relationship~\eqref{tau_equation} then follows from equation~\eqref{tau_formula}.

\section*{Acknowledgments}
JJHS was supported by the Maxwell Institute Graduate School in Analysis and its Applications, a Centre for Doctoral Training funded by the UK Engineering and Physical Sciences Research Council (grant EP/L016508/01), the Scottish Funding Council, Heriot-Watt University and the University of Edinburgh.

\bibliographystyle{plain}
\bibliography{Computing_extinction_path}

\end{document}